\documentstyle{amsppt}
\magnification=1200
\hsize=150truemm
\vsize=224.4truemm
\hoffset=4.8truemm
\voffset=12truemm

\NoRunningHeads

 \newsymbol\lessim 132E
\define\C{{\bold C}}
 
\define\R{{\bold R}}
\define\Z{{\bold Z}}
 
\let\thm\proclaim
\let\fthm\endproclaim
\define\bo{\partial  }  
 
 \define\p{$\C^2$}
 \define\s{$S^3 $}
\define\ba{$B^4 $}
\define\cb{$\Cal B$}
\define\cc{$\Cal C$}
\newcount\tagno
\newcount\secno
\newcount\subsecno
\newcount\stno
\global\subsecno=1
\global\tagno=0
\define\ntag{\global\advance\tagno by 1\tag{\the\tagno}}

\define\sta{\ 
{\the\secno}.\the\stno
\global\advance\stno by 1}

\define\stas{\the\stno
\global\advance\stno by 1}

\define\sect{\global\advance\secno by 1
\global\subsecno=1\global\stno=1\
{\the\secno}. }

\def\nom#1{\edef#1{{\the\secno}.\the\stno}}
\def\inom#1{\edef#1{\the\stno}}
\def\eqnom#1{\edef#1{(\the\tagno)}}

\newcount\refno
\global\refno=0

\def\nextref#1{
      \global\advance\refno by 1
      \xdef#1{\the\refno}}

\def\bref {\ref\global\advance\refno by 1\key{\the\refno}}


\nextref\A
\nextref\DI
\nextref\BK
\nextref\GO
\nextref\GR
\nextref\K
\nextref\KR
\nextref\SI
\nextref\ST

\topmatter
\title 
Riemann surfaces and totally real tori
 \endtitle
\author Julien Duval
 and Damien Gayet\footnote"*"{supported by the ANR\newline}
\endauthor
\abstract Given a generic totally real torus unknotted in the unit sphere \s\ of \p, we prove the following alternative : either there exists a filling of the torus by holomorphic discs and the torus is rationally convex, or its rational hull contains a holomorphic annulus.

\null \noindent

MSC 2010 : 32E20, 58J32

Keywords : totally real torus, filling by holomorphic discs, rational convexity
\endabstract 
  \endtopmatter 
\document
 \subhead 0. Introduction \endsubhead
\stno=1

\null

\noindent In this paper we address the following question :
\noindent
given a totally real torus in \p, does there always exist a compact Riemann
surface in \p\  with boundary in the torus ?

\null

Recall that (closed connected) surfaces in \p\ are {\it totally real} if they are never tangent to a complex line. The only orientable ones are tori. Special cases are {\it lagrangian} tori, those 
on which the standard K\"ahler form of \p\ vanishes. 

\null

Our question is motivated by geometric function theory. Given a compact set $K$ in \p, its {\it polynomial hull} $\hat K$ is defined as $\hat K=\{z \text { in } $\p / $\vert P(z)\vert \leq \Vert P \Vert_K $ for every polynomial $ P\}$. The set $K$ is {\it polynomially convex} if $\hat K=K$. In this case $K$ satisfies Runge theorem. Remark that any compact Riemann surface with boundary in $K$ is contained in $\hat K$. It is therefore tempting to explain the presence of a non trivial hull by Riemann surfaces, at least for nice sets like orientable surfaces (they are not polynomially convex for homological reasons). But quite often a complex tangency of a surface gives locally birth to small holomorphic discs with boundary on it. Thus the very first global problem arises with totally real orientable surfaces, namely tori.

Note that, in the definitions above, instead of polynomials we could as well work with rational functions without poles on $K$. This gives rise to the notions of 
{\it rational hull} and {\it rational convexity}. 

\null
Here is a bit of history around our question. In 1985 Gromov [\GR] gave a positive answer for lagrangian tori, constructing holomorphic discs with boundary in them. In 1996 by the same method Alexander [\A] exhibited for every totally real torus a proper holomorphic disc with all its boundary except
one point in the torus. Later on [\DI] he gave examples of totally real tori without holomorphic discs with full boundary in them, but still presenting holomorphic annuli.

\null
In the present work we focus\footnote"*"{following [\KR] which by the way seems uncorrect (see our example)} on tori in the unit sphere \s\ of \p. They are {\it unknotted} if they are isotopic to the standard torus in \s. We prove the following

\thm{Theorem} Let $T$ be a generic totally real torus unknotted in \s. Then either it bounds in \ba\ a solid torus foliated by holomorphic discs and $T$ is rationally convex, or its rational hull contains a holomorphic annulus with boundary in $T$.
\fthm
The solid torus is called a {\it filling} of $T$ which is said in this case {\it fillable}. The theorem applies to an open dense subset in the space of totally real
tori unknotted in \s. With some more work we could remove this genericity condition (in the second case a pair of holomorphic discs could also show up) but for simplicity we will stick to our statement.

The standard torus is an example of the first situation, while the second is illustrated by the following 

\null\noindent
{\bf Example} (compare with [\DI]). Consider the conjugate Hopf fibration 
$$\align
\pi: S^3\subset \C^2 &\to S^2\subset \C\times \R \\
(z,w) &\mapsto (2zw,\vert z\vert ^2-\vert w\vert^2). \endalign
$$
Remark that the fibers of
$\pi$ are circles.
Denote by $T_\gamma$ the preimage by $\pi$ of an embedded closed curve $\gamma$ in $S^2$. Then $T_\gamma$ is an unknotted torus in \s, totally real if the projection of $\gamma$ on $\C$ is immersed. Choose this projection as a figure eight avoiding the origin. It follows (see [\DI]) that every compact Riemann surface with boundary in $T_\gamma$ is in a fiber of the polynomial $p(z,w)=2zw$. But $T_\gamma$ does not separate $p^{-1}(a)$ except if $a$ is the double point of the figure eight. We then get only one holomorphic annulus with boundary in $T$.

\null

The proof of the theorem relies on the technique of filling spheres by holomorphic discs in its ultimate form due to Bedford and Klingenberg [\BK] (see also [\K]). This is where the restriction to \s\ enters. Here is the scheme of the argument.

\null

Take any unknotted torus $T$ in \s\. It divides \s\ in two solid tori. In the same manner its hull $\hat T$ separates the unit ball \ba\ in two pseudoconvex components. At least one of them has a universal covering which unwinds the corresponding solid torus. Push slightly $T$ in this good component (hence the genericity). We therefore get as a lifting of $T$ a periodic cylinder $\tilde T$ sitting in a pseudoconvex boundary. Approach this infinite cylinder by a sequence of spheres $S_n$ containing more and more periods of the cylinder. 

We are now in position to apply the technique of filling. It provides a sequence of balls bounded by $S_n$ and foliated by holomorphic discs, whose boundaries define a foliation on a bigger and bigger part of $\tilde T$. This sequence of foliations converges to a periodic foliation \cb. The alternative reads as follows : either all the leaves of \cb\ are closed or there are (non compact) periodic leaves. In the former case the closed leaves bound holomorphic discs which project down to the filling of $T$. In the later we prove that the periodic leaves bound periodic holomorphic strips which project down to holomorphic annuli.

\null 

Before entering the details of the proof, we collect some background. In the sequel all submanifolds are embedded
if not otherwise mentioned.

\null

\subhead 1. Background \endsubhead

\null\noindent
a) {\bf Filling spheres.} 

\null\noindent
Recall Bedford and Klingenberg theorem [\BK] (see also [\K]).

\thm{Theorem} Let $\Omega$ be a bounded strictly pseudoconvex domain in $\C^2$ and $S$ a sphere in $\bo \Omega$. Suppose that the complex tangencies of $S$ are elliptic or (good) hyperbolic points.
Then $S$ bounds a unique ball $\Sigma$ in $\Omega$ foliated by holomorphic discs.
\fthm

This ball is called the {\it filling} of $S$. The {\it complex tangencies} of $S$ are the points where $S$ is tangent to a complex line. Being of elliptic or (good) hyperbolic type (see [\BK] for the definition) is a generic condition.

\null

The picture looks as follows. Take a sphere in $\R^3$ endowed with its height function, which is Morse if the sphere is generic.
 Elliptic points correspond to local maxima and minima of the height,
while hyperbolic points translate in saddle points. The filling corresponds to the ball bounded by the sphere foliated by the level sets of the height.
Therefore all the holomorphic discs of the filling are smooth up to the boundary except those touching a hyperbolic point which present angles.

\null

Another way to describe the complex points of $S$ is via its {\it characteristic foliation} \cc. This is the foliation generated by the line field 
$T_{\C}\bo \Omega\cap TS$
where $T_{\C}\bo \Omega$ is the complex part of $T\bo \Omega$. It is singular precisely at the complex points of $S$, elliptic points corresponding
to foci and hyperbolic to saddle points.

\null

The foliation by holomorphic discs of the filling gives by restriction to its boundary another foliation of $S$, also singular at its complex points. These two foliations are transversal [\BK]. In particular the boundary of a holomorphic disc in the filling cuts at most once a leaf of \cc. 

\null
Here are further properties of the filling. First every compact Riemann surface in $\Omega$ with boundary in $S$ is contained in $\Sigma$. Next $\Sigma$ is the enveloppe of holomorphy of $S$. Hence $\Sigma$ is contained in any pseudoconvex domain containing $S$. Finally $\Sigma$ has a stability property. If we perturb slightly $S$ in $S'$ in $\bo \Omega$ then the filling of $S'$ is closed to $\Sigma$. Note that this stability also holds for the filling of a totally real torus in \s\ if it exists.

\null

In the sequel we will apply Bedford and Klingenberg technique to spheres in $\bo \tilde \Omega$ where $\tilde \Omega$ is the universal covering of a strictly pseudoconvex domain $\Omega$ in $\C^2$. The reader can check that all the arguments of [\BK] apply, mutatis mutandis.

\null\noindent
b) {\bf Holomorphic discs.} 

\null\noindent
We look more closely at holomorphic discs in $\tilde \Omega$.

We start with some length-area estimates.
Equip $\tilde \Omega$ with the lift of the standard K\"ahler form of $\C^2$. Note that it has a bounded primitive. Hence by Stokes theorem there is a constant $C$ such that for any holomorphic disc $\Delta$ in $\tilde \Omega$
$$ \text{area}(\Delta) \leq C\  \text{length}(\bo \Delta). \tag1$$

On the other hand the classical Ahlfors length-area estimate reads as follows.
Let $\Delta$ be a holomorphic disc in $\tilde \Omega$ parametrized by the upper half plane $H$ via a holomorphic map $f$ smooth up to $\R$. Denote by $a(r)$ the area of $f(H\cap D_r)$ and by $l(r)$
the length of $f(H\cap\bo D_r)$ where $D_r$ is the disc centered at $0$ of radius $r$. Then 
$$ \int_r^R\left(\frac {l}{a}\right)^2\frac{d\rho}{2\pi \rho} \leq \frac{1}{a(r)}. \tag2$$

As a consequence there exists a sequence $(r_n)$ of radii going to infinity such that $l(r_n)=o(a(r_n))$. If the area of the disc is finite we get then a sequence of arcs in its boundary which become asymptotically closed. Suppose the area infinite. According to $(1)$ we have $a(r_n)\leq C(l(r_n) + \lambda(r_n))$ where $\lambda(r)=$length$(f([-r,r]))$. We infer that $l(r_n)=o(\lambda(r_n))$. In other words we still get a sort of (infinite) {\it pinching} of the boundary : a sequence of arcs whose length goes to infinity faster than the distance between their extremities.

\null We turn now to sequences of holomorphic discs. Take a sequence parametrized
by $f_n : H \to \tilde \Omega$ where $f_n$ is holomorphic and smooth up to $\R$. Suppose moreover $f_n(i)$ converging. As $\Omega$ is bounded the derivative $\Vert f_n'\Vert$ remains locally uniformly bounded inside $H$. Hence $(f_n)$ converges (up to extraction) toward a holomorphic map $f: H\to \tilde \Omega$. 

\newpage
We focus on the behaviour of the sequence on $\R$ near the origin.
Suppose first that $a_n(\epsilon)$ blows up (we keep the notations above). By a similar argument there exists a sequence $(r_n)$ of radii between $\epsilon$ and $2\epsilon$ such that $l_n(r_n) = o(\lambda_n(r_n))$, hence again a pinching of the boundaries : a sequence of arcs $f_n([-r_n,r_n])$ whose length blows up faster than the distance between their extremities.

On the other hand if $a_n(\epsilon)$ is uniformly bounded and $f_n(\R)$ contained in a fixed totally real surface, then by Gromov compactness theorem [\SI] $f$ extends smoothly up to $]-\epsilon, \epsilon[$ with values in the surface. Moreover the convergence holds up to this interval except at finitely many points.

\null\noindent
c) {\bf Geometric function theory.} 

\null\noindent
We will use the following facts concerning polynomial convexity [\ST]. Let $K$ be a compact set in \s\ separating the sphere in finitely many components, then its polynomial hull $\hat K$ divides \ba\ in the same number of components. Moreover by Rossi local maximum principle these components are pseudoconvex domains.

We move on to rational convexity. The rational hull $r(K)$
of a compact set $K$ in $\C^2$ is geometrically defined as the set of points $z$ such that any algebraic curve passing through $z$ meets $K$. If $K \subset P$ where $P$ is a rational polyhedron, then the algebraic curves can be replaced by holomorphic curves in $P$. An obstruction to rational convexity is the presence of a compact Riemann surface with boundary in $K$ with the additional restriction that this boundary bounds in $K$. In our theorem (second situation) the holomorphic annulus will by construction satisfy this condition and therefore be part of $r(T)$. As for the first situation we have the following

\thm{Lemma} A fillable totally real torus in \s\ is rationally convex. \fthm

\noindent{\it Proof. } Call $T$ the torus and $\Sigma$ its filling. We first prove that $\Sigma$ is rationally convex. By Rossi local maximum principle and the Runge property of \ba\ it is enough to construct through any point near $\Sigma$ in the ball \ba\ a holomorphic curve in \ba\ (smooth up to \s) avoiding $\Sigma$. We produce them by stability of the filling. Foliate a neighborhood of $T$ in \s\ by tori, then the fillings of these tori foliate a neighborhood of $\Sigma$ in \ba. Therefore the corresponding holomorphic discs fill out this neighborhood and avoid $T$ if they are not in $\Sigma$. At this stage $r(T) \subset \Sigma$.

We prove now that $r(T)=T$. According to the first step $\Sigma$ is a decreasing limit of rational polyhedrons. It is then enough to construct through any point $z$ of $\Sigma \setminus T$ a holomorphic curve in a neighborhood of $\Sigma$ avoiding $T$. Take through $z$ a real closed curve in $\Sigma \setminus T$ transversal to the holomorphic discs, parametrized by the unit circle $\Gamma$. Extend this parametrization as a smooth map $f$ from a neighborhood of $\Gamma$ in such a way that $\bar \bo f$ vanishes to infinite order along $\Gamma$. By solving an adequate $\bar \bo$-equation perturb now $f$ into a holomorphic map. This map parametrizes a thin holomorphic annulus still passing through $z$ and intersecting $\Sigma$ near the initial curve, hence avoiding $T$.

\null

\subhead 2. Proof of the theorem \endsubhead

\null\noindent
a) {\bf The set up.}

\null\noindent
Consider an unknotted totally real torus $T$ in \s. It divides \s\ into two solid tori $\omega_i$ and its polynomial hull $\hat T$ 
separates \ba\ into two pseudoconvex domains $\Omega_i$ containing $\omega_i$ in their closure (\S1 c)). We have the following

\thm{Lemma} For one of these domains the map $H_1(\omega_i,\Z) \to 
H_1(\bar\Omega_i,\Z)$ is injective. \fthm
\noindent {\it Proof. }
If not, denote by $\gamma_i$ a generator of $H_1(\omega_i,\Z)$. By assumption a multiple of it bounds
in $\Omega_i$. For simplicity suppose that $\gamma_i$ itself is the boundary of a singular chain $C_i$ in
 $\Omega_i$. On the other hand $\gamma_1$ and $\gamma_2$ are linked in \s\ : one can find 
a disc $D_1$ in \s\ whose boundary is $\gamma_1$ and cutting $\gamma_2$ once. Take now a disc $D_2$ with boundary $\gamma_2$ outside the unit ball. We then get two cycles $C_i\cup D_i$ in $\C^2$ which
intersect only once. This is impossible.

\null

To fix the ideas let $\Omega_1$ be this good side. Being pseudoconvex it can be approximated from the inside by strictly pseudoconvex domains. Actually as $\Omega_1$ is already strictly convex along $\omega_1$ we don't need to modify it there. We get in this way a strictly pseudoconvex domain
$\Omega$ approximating $\Omega_1$ which contains in its boundary a solid torus $\omega$ slightly smaller than $\omega_1$. The injectivity of $H_1(\omega,\Z) \to 
H_1(\Omega,\Z)$ remains true.
 We push slightly the torus inside $\omega$, still calling it $T$.

\null

Consider the universal covering $p:\tilde \Omega\to \Omega$. Because $\pi_1(\omega) \to
\pi_1(\Omega)$ is injective, all the components of $p^{-1}(\omega)$ are diffeomorphic to $\R\times D^2$. Fix one of them and call it $\tilde \omega$.
Then $T$ lifts to a cylinder $\tilde T$ (diffeomorphic to $\R\times S^1$) inside $\tilde \omega$.
Let $\tau$ be the automorphism of $\tilde \Omega$ induced by the action of a generator of $\pi_1(\omega)$.
It acts on $\tilde \omega$ as a translation on the factor $\R$ and $\tilde T$ is invariant under this action.

\null
 
Construct the spheres $S_n$ approximating the infinite cylinder $\tilde T$ as follows. Consider a
 cylinder $C$ which is a fundamental
domain of the action of $\tau$ in $\tilde T$. Denote by $c$ and $\tau (c)$ its boundaries. Take a disc $\Delta$ bounding $c$ contained in the solid cylinder bounded by $\tilde T$ in $\tilde \omega$. Define the cap $\Delta^-$ (resp. $\Delta^+$) as the smoothing of the disc with corner $C\cup \Delta$ (resp. $\tau^{-1}(C)\cup \Delta$). Note that their complex points can be made generic after perturbation. The sphere $S_n$ is simply $(\cup_{i=-n}^{n-1}\ \tau^i(C) )\cup \tau^{n+1}(\Delta^+)\cup \tau^{-n-1}(\Delta^-)$. Call $\Sigma_n$ its filling.

\null\noindent
b) {\bf Convergence of the fillings.}

\null\noindent
 We want to investigate the limit of $\Sigma_n$. Let $B_n$ the ball bounded by $S_n$ in $\tilde \omega$, and $\Theta_n$ the bounded component of the complement of $\Sigma_n \cup B_n$ in $\tilde  \Omega$.

 We prove first that $\Theta_n$ increases toward a domain $\Theta$ invariant by $\tau$. Indeed deform slightly $S_{n+1}$ in $S_{n+1}'$ inside $\tilde \omega$ in such a way that $B_{n+1}' \supset S_n$. Denote by $\Sigma_{n+1}'$ the filling of
$S_{n+1}'$ and $\Theta_{n+1}'$ the corresponding bounded component. It is pseudoconvex and its boundary contains $S_n$. Hence $\Sigma_n \subset \Theta_{n+1}'$, so $\Theta_n\subset\Theta_{n+1}'$ and $\Theta_n\subset\Theta_{n+1}$ by deforming $S_{n+1}'$ back to $S_{n+1}$ (see \S1 a)). The invariance by $\tau$ of $\Theta$ is along the same lines, using the fact that $\tau^{\pm1}(S_{n-1}) \subset B_{n+1}'$.

\null

Denote by $\Cal B_n$ the foliation of $S_n$ given by the boundaries of the holomorphic discs in $\Sigma_n$.

\thm{Lemma} The foliations $\Cal B_n \vert_{\tilde T}$ converge toward a foliation \cb\ of $\tilde T$, invariant by $\tau$. Moreover the leaves of \cb\ extend slightly inside $\tilde \Omega$ as thin holomorphic strips.
\fthm
\noindent {\it Proof.} Given a point in $\tilde T$ denote by $l_n$ the tangent line of $\Cal B_n$ at this point for large $n$.
We know that the holomorphic discs whose boundaries give $\Cal B_{n+1}$ are outside $\Theta_n$.
This means that $l_n$ always rotates in the same direction when $n$ increases. On the other hand $l_n$ is transversal to the tangent of the characteristic foliation
of $\tilde T$ at this point. Hence $l_n$ reaches a limit position as $n$ goes to infinity. In this way we define a limit line field $L$ on $\tilde T$ which is invariant by $\tau$ (see above). 

Let us verify now that $L$ integrates in a foliation by curves bounding thin holomorphic strips inside $\tilde \Omega$. Locally on $\tilde T$ the curves of $\Cal B_n \vert_{\tilde T}$ are uniformly transversal to the characteristic foliation [\BK]. Therefore they are graphs
of uniformly Lipschitz functions which converge. Hence $L$ integrates
in Lipschitz curves forming the foliation \cb. Moreover this uniform transversality extends slightly in the interior of $\tilde \Omega$ [\BK]. Locally near $\tilde T$ the holomorphic discs of $\Sigma_n$ are graphs of holomorphic functions on domains in a fixed complex line which converge. Hence the thin holomorphic strips along the curves of \cb, and the regularity of these curves by ellipticity of the $\bar \bo$-operator with totally real boundary conditions [\SI]. Note that \cb\ is only transversally continuous.

\null\noindent
c) {\bf The alternative.}  

\null \noindent
Consider in $\Cal B_n$ the curve $\delta_n$ passing through a given point of $c$ (the curve dividing out $S_n$ in two equal parts). 

The alternative reads as follows : either the length of $\delta_n$ does not blow up or it does.

We will see that it can also be worded in terms of \cb\ :
either all its leaves are closed or it presents periodic leaves (invariant by $\tau$).

\null
The easiest case of the alternative is the first. Because its length remains bounded, $\delta_n$ does not visit the caps of $S_n$ for large $n$. Hence we get a holomorphic disc $\Delta$ of $\Sigma_n$ whose boundary remains in $\tilde T$. Now take $p$ such that $\Delta$ and $\tau^p(\Delta)$ are disjoint. Applying Bedford and Klingenberg technique [\BK] we can interpolate between them by a family of holomorphic discs in $\tilde \Omega$ with boundaries in $\tilde T$ forming a solid cylinder. This solid cylinder is part of $\Sigma_m$ for large $m$ (\S1 a)). In other words $\Cal B_m$ stabilizes between $\bo \Delta$ and $\tau^p(\bo \Delta)$. The foliation \cb\ is then nothing else than the foliation by the boundaries of this family of holomorphic discs, extended to $\tilde T$ by the action of $\tau^p$. The family is actually invariant by $\tau$ because \cb\ is, and we get a solid cylinder foliated by holomorphic discs between $\Delta$ and $\tau(\Delta)$ which projects down to the filling of $T$.

\null
The second case is more delicate. 

\null
We prove first that if the length of $\delta_n$ blows up then \cb\ has non compact proper leaves and even periodic ones.

This requires an extra discussion of the characteristic foliation of $T$. By Denjoy theorem [\GO] any smooth foliation on $T$ can be perturbed in order to get only a finite number of attracting or repulsive cycles (closed leaves). As we want to analyze boundaries in $T$ of holomorphic discs in $\Omega$ which are uniformly transversal to the characteristic foliation, we may as well replace the characteristic foliation by such a perturbation, which we do henceforth. 

\null
Note that the characteristic cycles do not lift as closed curves to $\tilde T$. If it were the case the spheres $S_n$ (for large $n$) would get a characteristic cycle. But the foliation $\Cal B_n$ given by the filling can be described as the levels of a function $\phi$ (\S1 a)). Necessarily $\phi$ would have a critical point on the cycle, contradicting the transversality between the two foliations. Hence the characteristic foliation \cc\ on $\tilde T$ has finitely many periodic attracting or repulsive leaves. 
Remove thin tubes along these periodic leaves from $\tilde T$. We get a finite number of strips which can be parametrized by $\R \times[0,1]$ via a periodic diffeomorphism sending the vertical foliation to the characteristic foliation. As a consequence we show that $\delta_n$ cannot remain in a compact part of $\tilde T$. Indeed we know that $\delta_n$ crosses at most once each periodic leaf (\S1 a)). By uniform transversality this implies that $\delta_n$ cuts each thin tube in at most one short arc. Hence $\delta_n$ spends the major part of its length in the strips where it is uniformally transversal to the vertical foliation (via the parametrization). Therefore it cannot remain in a compact part of $\tilde T$. Passing to the limit we get at least one non compact proper leaf $\gamma$ in \cb\ whose ends are contained in strips. 

By the way a similar discussion shows that a pinching of boundaries in $\tilde T$ (in the sense of \S1 b)) has to cross at least one periodic leaf of \cc.

\null

We prove now the existence of periodic leaves in \cb. For this take one end of $\gamma$ in the parametrization of its strip. Being transversal to the vertical foliation it is a graph of a map $h: [a,+\infty[\to [0,1]$. As $\gamma$ is disjoint from $\tau(\gamma)$ we get that $h\circ t^{-1}$ is say larger than $h$ on $[a,+\infty[$ if $\tau$ gives a positive translation $t$ in the parametrization.
This means that $h\circ t^{-n}$ is increasing toward a limit $g$ invariant by $t$.
Its graph gives a periodic leaf $\alpha$ of \cb\.

 \null

To end the proof we glue to $\alpha$ a global holomorphic strip. For this we go back to the construction of \cb. Recall that $\alpha$ is the limit of a part of $\bo \Delta_n$ where $\Delta_n$ is a holomorphic disc in the filling $\Sigma_n$. Moreover $\alpha$ bounds already a thin holomorphic strip which is the limit of a part of $\Delta_n$. Fix two points $x$ on $\alpha$ and $y$ in the interior of this strip. They are limits of $x_n$ on $\bo \Delta_n$ and $y_n$ in the interior of $\Delta_n$. Parametrize $\Delta_n$ by a holomorphic map $f_n :(D,\bo D) \to (\tilde \Omega, S_n)$, where $D$ is the unit disc, $f_n$ smooth up to $\bo D$ and $f_n(0)=y_n,f_n(1)=x_n$.

\thm{Lemma} The sequence $(f_n)$ converges to a smooth map $f :(D, \bo D\setminus$ two points$) \to(\tilde \Omega, \tilde T)$ whose image $S$ is invariant by $\tau$. \fthm

Hence $S$ projects down in $\Omega$ to a holomorphic annulus $A$ such that $\bo A\subset T$ and $\bo A$ bounds an annulus in $T$, which is what we wanted.

\null
\noindent{\it Proof. }  We begin by proving that $(f_n)$ converges up to $\bo D$ in a neighborhood of $1$. We parametrize the part of $\Delta_n$ converging to a given compact piece of the thin holomorphic strip along $\alpha$ by the upper half unit disc $D^+$ via $h_n$ in such a way that $h_n(0)=x_n$ and $h_n^{-1}(y_n)$ has a limit in the interior of $D^+$. By Schwarz reflection $f_n^{-1}\circ h_n : (D^+,]-1,1[)\to (D,\bo D)$ converges toward a non constant map $k$. Hence $f_n$ converges on the image of $k$. Taking an increasing sequence of compact pieces to the whole thin strip we get that $(f_n)$ converges to $f$ on a neighborhood of an arc $a\subset \bo D$ such that $f\vert_a$ parametrizes $\alpha$. Extend $f$ to $D$ (\S1 b)) and write $S$ for $f(D)$.

\null
Next we prove that the complement of $a$ in $\bo D$ is not reduced to a single point. If it were the case we could parametrize $S$ by $H$ instead of $D$ putting this point at infinity and get (\S1 b)) a pinching of $\alpha$, a sequence of subarcs $\alpha_n$ whose length blows up faster than the distance between their extremities. But this is impossible because this distance is comparable to the number of periods in $\alpha_n$, hence to its length. 

\null
At this stage we have already the convergence of $(f_n)$ toward $f$ on $D$ up to the arc $a$. It remains to see what happens on $b=\bo D\setminus \bar a$. Recall (\S1 b)) that to extend $f$ to $b$ it is enough to ensure an area bound for $f_n$ near this arc. 

\null
We prove now that the area blows up near at most finitely many points of $b$, meaning that $f$ extends except at these points. For this recall that each blow-up gives rise to a pinching of $\bo \Delta_n$ (\S1 b)). Recall also that $\bo \Delta_n$ cuts at most once each leaf of \cc\ (\S1 a)). We have two possibilities for a pinching.
Either it remains in $\tilde T$ or it visits a cap of $S_n$ for large $n$. In the former case it has to cross a periodic leaf of \cc. There are finitely many of them. In the later it has to cross the boundary of the cap. Again this bounds a priori their number. If not, more and more strands of $\bo \Delta_n$ would accumulate somewhere at the boundaries of the caps, meaning that $\bo \Delta_n$ would cut more than once the characteristic leaf through this point. 

\null
Lastly note that $S$
is invariant by $\tau$ because $\alpha$ is.
Hence $\tau$ translates via $f$ in an automorphism of $D$ globally preserving $a$ and $b$, a hyperbolic translation which has infinite orbits. This implies that $f$ actually extends to the whole $b$, which concludes.
 
\null\noindent{\bf Final remark.} This method allows us to glue to each leaf of \cb\ a holomorphic disc with finitely many singular points at the boundary. Their union projects down to a set looking very much like a Levi-flat hypersurface bounding $T$, except for its regularity.

\Refs
\widestnumber\no{99}
\refno=0

\bref \by H. Alexander\paper Gromov's method and Bennequin's problem\jour Inv. Math.\vol125\yr1996\pages135--148
\endref

\bref \by H. Alexander\paper Disks with boundaries in totally real and Lagrangian manifolds\jour Duke Math. J.\vol100\yr1999\pages131--138
\endref

\bref \by E. Bedford and W. Klingenberg\paper On the envelope of holomorphy of a 2-sphere in $\C^2$\jour J. Amer. Math. Soc.\vol4\yr1991\pages623--646
\endref

\bref \by C. Godbillon\book Dynamical systems on surfaces \bookinfo Universitext\publ Springer \yr1983 \publaddr Berlin
\endref

\bref \by M. Gromov\paper Pseudoholomorphic curves in symplectic manifolds\jour Inv. Math.\vol82\yr1985\pages307--347
\endref

\bref \by N. G. Kruzhilin\paper Two-dimensional spheres on the boundaries of pseudoconvex domains in $\C^2$\jour Math. USSR-Izv.\vol39\yr1992\pages1151--1187
\endref

\bref \by N. G. Kruzhilin\paper Holomorphic disks with boundaries in totally real tori in $\C^2$\jour Math. Notes\vol56\yr1994\pages1244--1248
\endref

\bref \by J.-C. Sikorav \book Some properties of holomorphic curves in almost complex manifolds, {\rm in} Holomorphic curves in symplectic geometry
\pages165--189 \bookinfo Prog. Math. \vol117 \publ Birkh\"auser \yr1994 \publaddr Basel
\endref

\bref \by E. L. Stout \book Polynomial convexity \bookinfo Prog. Math. \vol261 \publ Birkh\"auser \yr2007 \publaddr Boston
\endref

\endRefs
\address J. Duval, Laboratoire de math\'ematique, Universit\'e Paris-Sud, B\^atiment 425, 91405 Orsay Cedex, France\endaddress

\email julien.duval\@math.u-psud.fr \endemail

\address D. Gayet, Institut Camille Jordan, CNRS UMR 5208, UFR de math\'ematiques, Universit\'e Lyon I, B\^atiment Braconnier, 69622 Villeurbanne Cedex, France\endaddress

\email gayet\@math.univ-lyon1.fr\endemail
\enddocument